\documentclass[journal,twoside,web]{ieeecolor}
\usepackage{generic}
\usepackage{cite}
\usepackage{amsmath,amssymb,amsfonts}
\usepackage{graphicx}
\usepackage{hyperref}
\usepackage{textcomp}
\usepackage{algpseudocode}
\usepackage{array}
\usepackage{bbm}
\usepackage{balance}
\usepackage{booktabs}
\usepackage{cases}
\usepackage{dsfont}
\usepackage{mathrsfs}
\usepackage{mathtools}
\usepackage{multirow}
\usepackage{makecell}
\usepackage{textcomp}
\usepackage{soul}
\usepackage{url}
\usepackage{siunitx}
\usepackage{newtxtext,newtxmath}

\newtheorem{proposition}{Proposition}

\newtheorem{theorem}{Theorem}
\newtheorem{defi}{Definition}
\newtheorem{ass}{Assumption}
\newtheorem{remark}{Remark}

\def\BibTeX{{\rm B\kern-.05em{\sc i\kern-.025em b}\kern-.08em
    T\kern-.1667em\lower.7ex\hbox{E}\kern-.125emX}}
\begin{document}
\title{\LARGE On Min-Max Robust Data-Driven Predictive Control Considering Non-Unique Solutions to Behavioral Representation}
\author{Yibo Wang, Qingyuan Liu, Chao Shang, \IEEEmembership{Member, IEEE}
    \thanks{This work was supported by National Natural Science Foundation of China under Grants 62373211 and 62327807. \textit{(Corresponding author: Chao Shang.)}}
    \thanks{Y. Wang, Q. Liu, and C. Shang are with Department of Automation, Beijing National Research Center for Information Science and Technology, Tsinghua University, Beijing 100084, China (e-mail: wyb21@mails.tsinghua.edu.cn, lqy20@mails.tsinghua.edu.cn, c-shang@tsinghua.edu.cn). }}

\maketitle

\begin{abstract}
    Direct data-driven control methods are known to be vulnerable to uncertainty in stochastic systems. In this paper, we propose a new robust data-driven predictive control (DDPC) framework. By analyzing non-unique solutions to behavioral representation, \textcolor{black}{we gain insight into the inherent lack of robustness in subspace predictive control (SPC) and its projection-based regularized variant.} This stimulates us to construct an uncertainty set that captures all admissible output trajectories deviating from nominal subspace predictions, which results in a min-max robust formulation of DDPC that endows control sequences with robustness against such unknown deviations. We establish theoretical performance guarantees under bounded additive noise and develop tractable convex reformulations. To mitigate the conservatism of robust design, a feedback robust DDPC scheme is further proposed by incorporating an affine feedback policy. Simulation studies show that the proposed methods effectively robustify SPC and outperform the projection-based regularization.
\end{abstract}

\begin{IEEEkeywords}
    Data-driven control, robust control, robust optimization.
\end{IEEEkeywords}

\section{Introduction}
\label{sec:introduction}
\IEEEPARstart{I}{n} recent years, data-driven control methods have received increasing attention for deriving the control policies directly from raw data while bypassing the cumbersome system identification \cite{shang2019data,markovsky2021behavioral}. Among these methods, a mainstream is built upon the well-known Willems' fundamental lemma \cite{willems2005note}, where all admissible trajectories of a deterministic linear time-invariant (LTI) system can be characterized by a single input-output sequence \cite{van2020willems}. This forms the foundation of the generic data-driven predictive control (DDPC) method, where a data-driven predictor based on a set of data equations is used for output predictions \cite{coulson2019data}.

Note that the generic DDPC method is only applicable to deterministic LTI systems. In the presence of uncertainty, the fundamental lemma is no longer strictly valid, \textcolor{black}{compromising control performance} \cite{yin2023maximum}. To attenuate the effect of uncertainty, several remedies based on generic DDPC have been proposed. Adopting the minimum-norm solution to the data equations yields the well-known Subspace Predictive Control (SPC) method \cite{favoreel1999spc}. \textcolor{black}{Alternatively, regularization terms such as the norm-based penalty \cite{coulson2019data,berberich2021data} or the projection-based regularizer \cite{dorfler2023bridging} have been introduced to balance control cost and prediction fitting. The innovation-based DDPC method has also been developed by characterizing the additive stochastic noise as innovation estimates \cite{wang2025data}.}

In this paper, we first analyze how uncertainty is addressed by some existing DDPC formulations from a new perspective, \textcolor{black}{that is the \textit{non-uniqueness} of solutions to behavioral representations}. This sheds some light on the inherent limitations of SPC and PBR-DDPC in addressing the uncertainty in dynamical systems. Specifically, the SPC uniquely chooses the minimum-norm solution, whereas the projection-based regularized DDPC (PBR-DDPC) only seeks the best-case control performance among a continuum of non-unique solutions in an optimistic sense. This renders PBR-DDPC sensitive to uncertainty and leads to a lack of robustness \cite{gotoh2025technical}.

Building upon these analyses, we propose a novel min-max robust formulation of DDPC to tackle uncertainty in dynamic systems. \textcolor{black}{Specifically, we construct an uncertainty set from non-unique solutions to behavioral representations,} which encloses a continuum of admissible output trajectories deviating from the ``nominal'' SPC predictions to some degree. This motivates a new min-max robust formulation of DDPC that optimizes the worst-case performance under all possible realizations in the constructed uncertainty set. Furthermore, under bounded additive input-output noise, we establish theoretical performance guarantees for the robust DDPC by appropriately sizing the uncertainty set, and derive a tractable convex semi-definite programming (SDP) problem corresponding to the min-max robust formulation. Additionally, to alleviate the conservatism of robust design, we propose a feedback robust DDPC scheme by incorporating an affine feedback control policy, which can appropriately describe the adaptation of the control policy to future uncertainty in receding horizon implementation. We also offer a performance guarantee for the feedback robust DDPC under bounded additive noise. A numerical example showcases that our methods achieve significant performance improvement over SPC, while PBR-DDPC cannot.

The rest of this paper is organized as follows. Section II gives a brief review of DDPC methods in noise-free and noisy cases. Section III presents the robust DDPC method. Section IV derives the feedback robust DDPC scheme with an affine feedback control policy. Section V provides the numerical studies. Finally, Section VI concludes the paper.

\textit{Notation:} Let $[a]$ denote the index set $\{1,2,\cdots,a\}$. For a sequence $\{x(i)\}_{i=1}^N$, $x_{[i:j]}$ denotes its restriction to the interval $[i,j]$ as ${\rm col}(x(i),\cdots, x(j))=\begin{bmatrix}x(i)^\top&\cdots&x(j)^\top\end{bmatrix}^\top$. This notation extends to sequences of matrices with identical column dimensions via vertical stacking. The operator $\mathcal{H}_s(x_{[i:j]})$ builds a block Hankel matrix of depth $s$ from the sequence. We use $\sigma_{i}(\cdot)$ and $\sigma_{\max}(\cdot)$ to denote the $i$-th and maximal singular values. \textcolor{black}{We use ${\rm im}(X)$ and $\mathcal{N}(X)$ for the image and right null spaces of matrix $X$, respectively,} and $\star$ to represent symmetric blocks in linear matrix inequalities (LMIs).

\section{Preliminaries}

\subsection{DDPC of Deterministic LTI Systems}

Consider a deterministic discrete-time LTI system with a minimal representation:
\begin{equation}
    \label{equation: innovation-form LTI system}
    \left \{
    \begin{aligned}
        \bar{x}(t+1) & =A\bar{x}(t)+B\bar{u}(t),  \\
        \bar{y}(t)   & =C\bar{x}(t)+D\bar{u}(t),
    \end{aligned} \right .
\end{equation}
where $\bar{x}(t)\in\mathbb{R}^{n_x}$, $\bar{u}(t)\in\mathbb{R}^{n_u}$, and $\bar{y}(t)\in\mathbb{R}^{n_y}$ denote state, input, and output, respectively. Throughout this paper, we use the topmark $\bar{\cdot}$ to denote data vectors and matrices based on noiseless data from \eqref{equation: innovation-form LTI system}. \textcolor{black}{Let an input sequence $\{\bar{u}(i)\}_{i=1}^N$ be persistently exciting of order $L+n_x$ (i.e., its Hankel matrix of depth $L+n_x$ has full row rank \cite{willems2005note}).} By the fundamental lemma, any $L$-long input-output trajectory $\{\bar{u}(i),\bar{y}(i)\}_{i=t-L_p}^{t+L_f-1}$ with $L=L_p+L_f$ is a trajectory of system \eqref{equation: innovation-form LTI system} if and only if there exists \textcolor{black}{$g\in\mathbb{R}^{T}$ with $T=N-L+1$}, such that \cite{van2020willems}:
\begin{equation}
    \label{equation: data-based equations noise free}
    \begin{bmatrix}
        \mathcal{H}_L(\bar{u}^d_{[1:N]}) \\
        \mathcal{H}_L(\bar{y}^d_{[1:N]})
    \end{bmatrix}g=
    \begin{bmatrix}
        \bar{U}_d \\\bar{Y}_d
    \end{bmatrix}g=
    \begin{bmatrix}
        \bar{u}_{[t-L_p:t+L_f-1]} \\\bar{y}_{[t-L_p:t+L_f-1]}
    \end{bmatrix}.
\end{equation}
The relation \eqref{equation: data-based equations noise free} offers a \textit{behavioral data-driven} representation of the LTI system without knowing system matrices $\{A,B,C,D\}$ and paves the way for model-free output prediction. If $L_p\ge\ell(A,B,C,D)$ and $L_f\ge1$, \textcolor{black}{where $\ell$ is the lag of the system,} the future output $\bar{y}_f=\bar{y}_{[t:t+L_f-1]}$ at time instance $t$ can be uniquely determined by $\bar{y}_f=\bar{Y}_f g$ provided that $g$ solves the following linear equations:
\begin{equation}
    \label{equation: calculate g noise free}
    \underbrace{{\rm col}(
            \bar{U}_p,\bar{U}_f,\bar{Y}_p
        )}_{\triangleq\bar{\Phi}}g={\rm col}(
        \bar{u}_p,\bar{u}_f,\bar{y}_p),
\end{equation}
with
\begin{equation}
    \begin{aligned}
        &\bar{U}_p=\mathcal{H}_{L_p}(\bar{u}^d_{[1:N-L_f]}),~\bar{U}_f=\mathcal{H}_{L_f}(\bar{u}^d_{[L_p+1:N]}),\\
        &\bar{u}_p=\bar{u}_{[t-L_p:t-1]},~\bar{u}_f=\bar{u}_{[t:t+L_f-1]},
    \end{aligned}
\end{equation}
and $\{\bar{Y}_p,\bar{Y}_f,\bar{y}_p\}$ defined similarly \cite{markovsky2008data}. Using $\bar{y}_f=\bar{Y}_fg$ and \eqref{equation: calculate g noise free} as the data-driven output predictor, the generic DDPC can be framed as the following optimization problem \cite{coulson2019data}:
\begin{equation}
    \label{equation: generic DeePC}
    \begin{aligned}
        \min_{\bar{u}_f,\bar{y}_f,g} & ~\mathcal{J}(\bar{u}_f,\bar{y}_f)                                                              \\
        {\rm s.t.}~                  & \bar{\Phi}g=\bar{\phi}(\bar{u}_f),~\bar{y}_f=\bar{Y}_fg, \\
                                     & \bar{u}_f\in\mathbb{U},~\bar{y}_f\in\mathbb{Y},
    \end{aligned}
\end{equation}
where $g$ is a decision variable to be optimized together with $\{\bar{u}_f,\bar{y}_f\}$, \textcolor{black}{and $\bar{\phi}(\bar{u}_f)\triangleq{\rm col}(\bar{u}_p,\bar{u}_f,\bar{y}_p)$ is defined as the stacked data vector depending on the future input $\bar{u}_f$.} Here, we assume that the future inputs and outputs are confined to intersections of ellipsoids as $\mathbb{U}=\{u_f~|~\|G_u^i u_f+c_u^i\|^2\le1,~\forall i\in[l_u]\}$ and $\mathbb{Y}=\{y_f~|~\|G_y^i y_f+c_y^i\|^2\le1,~\forall i\in[l_y]]\}$. The cost function in \eqref{equation: generic DeePC} can be chosen as $\mathcal{J}(\bar{u}_f,\bar{y}_f)=\sum_{i=0}^{L_f-1}\|\bar{y}(t+i)-y_r(t+i)\|_Q^2+\|\bar{u}(t+i)\|_R^2$ with $Q,R\succ0$, where $y_r(\cdot)$ is the output reference.

% \begin{equation}
%     \label{equation: constraints on input and output}
%     \begin{aligned}
%         \mathbb{U} & =\{u_f~|~\|G_u^i u_f+c_u^i\|^2\le1,~\forall i\in[l_u]\},              \\
%         \mathbb{Y} & =\{y_f~|~\|G_y^i y_f+c_y^i\|^2\le1,~\forall i\in[l_y]]\},
%     \end{aligned}
% \end{equation}
% where $G_u^i\in\mathbb{R}^{m_y^i\times n_uL_f}$, $G_y^i\in\mathbb{R}^{m_u^i\times n_yL_f}$, and $l_u,~l_y \ge 1$ stand for the numbers of constraints on future input and output, respectively. The cost function in \eqref{equation: generic DeePC} can be chosen as:
% \begin{equation}
%     \label{equation: cost function}
%     \mathcal{J}(\bar{u}_f,\bar{y}_f)=\sum_{i=0}^{L_f-1}\|\bar{y}(t+i)-y_r(t+i)\|_Q^2+\|\bar{u}(t+i)\|_R^2,
% \end{equation}
% where $y_r(\cdot)$ is the output reference, and $Q,R\succ0$ denote the weighting matrices.

\subsection{SPC and Regularized DDPC}
In practice, the presence of noise and disturbance leads to high variance in output prediction and renders the solution to \eqref{equation: generic DeePC} constructed from noisy data unreliable. Hereafter, we use the unmarked data vector and matrices to represent those constructed from noisy data. To mitigate the effect of uncertainty, the SPC adopts the minimum-norm solution $g_{\rm pinv}=\Phi^\dagger{\rm col}(u_p,u_f,y_p)$, yielding the SPC predictor $\hat{y}_f=Y_f\Phi^\dagger{\rm col}(u_p,u_f,y_p)$, and solves the following quadratic programming (QP) problem \cite{favoreel1999spc}:
\begin{equation}
    \label{equation: SPC}
    \begin{aligned}
        \min_{u_f,\hat{y}_f} & ~\mathcal{J}(u_f,\hat{y}_f)               \\
        {\rm s.t.}~          & \hat{y}_f=Y_f \Phi^\dagger\phi(u_f),      \\
                             & u_f\in\mathbb{U},~\hat{y}_f\in\mathbb{Y},
    \end{aligned}
\end{equation}
\textcolor{black}{where $\phi(u_f)={\rm col}(u_p,u_f,y_p)$, analogous to $\bar{\phi}(\cdot)$, denotes the stacked vector depending on $u_f$.}

Recently, various regularized extensions of DDPC have been put forward to attenuate the effect of uncertainty. A representative is the projection-based regularized DDPC (PBR-DDPC) \cite{dorfler2023bridging}, which yields the following regularized optimization problem:
\begin{equation}
    \label{equation: original regularized DDPC}
    \begin{aligned}
        \min_{u_f,\hat{y}_f,g} & ~\mathcal{J}(u_f,\hat{y}_f)+\lambda\cdot\|\Phi^\perp g\|_2^2 \\
        {\rm s.t.}~            & {\rm col}(\Phi,Y_f)g={\rm col}(\phi(u_f),\hat{y}_f),\\
                               & u_f\in\mathbb{U},~\hat{y}_f\in\mathbb{Y},
    \end{aligned}
\end{equation}
where $\lambda>0$ is a regularization parameter, and $\Phi^\bot\triangleq I-\Phi^\dagger\Phi$ is an orthogonal projector onto $\mathcal{N}(\Phi)$. Interestingly, \eqref{equation: original regularized DDPC} turns out to be a convex relaxation of the bilevel problem, where a multi-step predictor is first identified from data and then used for predictive control design \cite{dorfler2023bridging}. By varying $\lambda$, a flexible balance can be made between control cost minimization and fitting an output predictor from data. When $\lambda \to \infty$, \eqref{equation: original regularized DDPC} reduces to the design problem \eqref{equation: SPC} of SPC.

\section{A New Scheme for Robustifying DDPC}
\label{Sec: Robust DDPC}

In this section, we first shed light on the deficiency of SPC \eqref{equation: SPC} and PBR-DDPC \eqref{equation: original regularized DDPC} from a new perspective, and then propose a min-max robust formulation of DDPC to mitigate the effect of uncertainty.

\subsection{Non-Unique Solutions to Behavioral Representation}

We highlight that the behavioral representation \eqref{equation: calculate g noise free} in the noise-free case is essentially an \textit{under-determined} equation, whose solution is not unique. Mathematically, any admissible solution $g$ to \eqref{equation: calculate g noise free} can be parameterized as the sum of the minimum-norm solution and an orthogonal term \cite{mattsson2023regularization}:
\begin{equation}
    \label{equation: g full solution noise free}
    g=\bar{\Phi}^\dagger\bar{\phi}(\bar{u}_f)+\bar{\Phi}^\perp w,
\end{equation}
where $\bar{\Phi}^\perp=I-\bar{\Phi}^\dagger\bar{\Phi}$, and $w\in\mathbb{R}^{T}$ is a free variable. In \eqref{equation: g full solution noise free}, the term $\bar{\Phi}^\dagger{\rm col}(\bar{u}_p,\bar{u}_f,\bar{y}_p)\in{\rm im}(\bar{\Phi}^\top)$, i.e., the row space of $\bar{\Phi}$, and the term $\bar{\Phi}^\bot w\in\mathcal{N}(\bar{\Phi})$, i.e., the right null space of $\bar{\Phi}$. Note that any choice of $w$ in \eqref{equation: g full solution noise free} ensures the satisfaction of \eqref{equation: calculate g noise free}. Put differently, one can by no means identify the value of $w$ from data, and this leads to the non-uniqueness of $g$ in $\mathcal{N}(\bar{\Phi})$. Despite this, the future output can be uniquely determined whenever exact data are collected from \eqref{equation: innovation-form LTI system} and satisfy the assumptions of the fundamental lemma:
\begin{equation}
    \label{equation: OP noise free}
    \begin{aligned}
        \bar{y}_f=\bar{Y}_f\bar{\Phi}^\dagger\bar{\phi}(\bar{u}_f)+\underbrace{\bar{Y}_f\bar{\Phi}^\bot}_{=0}w= \bar{Y}_f\bar{\Phi}^\dagger\bar{\phi}(\bar{u}_f),
    \end{aligned}
\end{equation}
where $\bar{Y}_f\bar{\Phi}^\bot = 0$ is a natural implication from the fundamental lemma. Indeed, this is also implied by the subspace identification literature since $\bar{Y}_f$ can always be expressed as a linear combination of $\{\bar{U}_p,\bar{U}_f,\bar{Y}_p\}$ \cite{huang2008dynamic} and thus $\bar{Y}_f\bar{\Phi}^\perp=0$. However, in practical dynamical systems subject to uncertainty, it is always the case that $Y_f\Phi^\bot\neq0$ and the fundamental lemma no longer holds strictly. In this case, the expression of $g$ reads as:
\begin{equation}
    \label{eq: 13}
    g=\Phi^\dagger\phi(u_f)+\Phi^\bot w,
\end{equation}
where $\Phi^\bot\Phi^\dagger=0$. Then the data-driven output predictor in \eqref{equation: OP noise free} takes the form:
\textcolor{black}{
\begin{equation}
    \label{equation: OP noise}
    \hat{y}_f(u_f,w)=Y_fg=Y_f\Phi^\dagger\phi(u_f)+Y_f\Phi^\bot w,
\end{equation}
where the output prediction is expressed as a function of $u_f$ and $w$, consisting of the nominal SPC predictor plus an additional term in $w$}. In this sense, the non-uniqueness of $g$ in the nullspace of $\Phi$ brings some uncertainty to the output prediction $\hat{y}_f$.

From this point of view, we revisit the ideas underlying SPC \cite{favoreel1999spc} and PBR-DDPC \cite{dorfler2023bridging}. In essence, the SPC predictor in \eqref{equation: SPC} simply assumes $w=0$ and that $g$ has no components in $\mathcal{N}(\Phi)$, which yields a unique pseudo-inverse solution $\Phi^\dagger\phi(u_f)$. As for the regularized problem \eqref{equation: original regularized DDPC} of PBR-DDPC, it tackles the non-uniqueness of $g$ in \eqref{eq: 13} by regarding $w$ as a free decision variable and consequently, allows the output prediction $\hat{y}_f$ to differ from the predictions from the SPC predictor in an \textit{optimistic} sense. Since $\Phi^\bot$ is idempotent, i.e., $\Phi^\bot \Phi^\bot = \Phi^\bot$, we can utilize the expression of $g$ in \eqref{eq: 13} to derive the equation $\Phi^\bot g=\Phi^\bot\left(\Phi^\dagger\phi(u_f)+\Phi^\bot w\right)=\Phi^\bot w$. This enables us to rewrite the regularized problem \eqref{equation: original regularized DDPC} as:
\begin{equation}
    \label{equation: regularized DDPC}
    \begin{aligned}
        \min_{u_f\in\mathbb{U},w} & ~\mathcal{J}\left(u_f,\textcolor{black}{\hat{y}_f(u_f,w)}\right)+\lambda\cdot\|\Phi^\bot w\|^2      \\
        {\rm s.t.}  ~~&~\textcolor{black}{\hat{y}_f(u_f,w)}\in\mathbb{Y},
    \end{aligned}
\end{equation}
where the function of output prediction $\hat{y}_f(u_f,w)$ is defined as \eqref{equation: OP noise}, and $w$ becomes a decision variable to be optimized. Using the duality of convex programming, one obtains the equivalence between \eqref{equation: regularized DDPC} and the following min-min problem:
\textcolor{black}{
\begin{equation}
    \label{equation: regularized DDPC with constraint}
    \begin{aligned}
        \min_{u_f\in\mathbb{U}}~\min_{w\in\mathcal{W}(\Lambda)}&~\mathcal{J}(u_f,\hat{y}_f(u_f,w))\\
        {\rm s.t.}~~&~\hat{y}_f(u_f,w)\in\mathbb{Y},
    \end{aligned}
\end{equation}
}
where the regularization term in \eqref{equation: regularized DDPC} now turns into an \textit{uncertainty set} of $w$ as $\mathcal{W}(\Lambda)\triangleq\{~w~|~\|\Phi^\bot w\|_2^2\le\Lambda~\}$. Here $\Lambda>0$ is the size parameter of $\mathcal{W}(\Lambda)$, which can also be interpreted as the \textit{maximal admissible deviation from the predictions of SPC predictor in \eqref{equation: SPC}}. Indeed, the regularization parameter $\lambda$ in \eqref{equation: regularized DDPC} can be interpreted as the optimal Lagrangian multiplier associated with the ellipsoidal constraint $w\in\mathcal{W}(\Lambda)$ in \eqref{equation: regularized DDPC with constraint}. For any finite $\lambda$ in \eqref{equation: regularized DDPC}, there always exists a corresponding $\Lambda$ such that optimal solutions to \eqref{equation: regularized DDPC with constraint} and \eqref{equation: regularized DDPC} coincide.

The formulation \eqref{equation: regularized DDPC with constraint} gives some new insights into the design problem \eqref{equation: original regularized DDPC} of PBR-DDPC. That is, solving \eqref{equation: original regularized DDPC} amounts to finding an input sequence $u_f^*$ along with a particular $w^* \in \mathcal{W}(\Lambda)$ that yields the \textit{lowest} control cost and fulfills all prescribed constraints. \textcolor{black}{This is essentially a criterion of \textit{optimism} since the selected input-output relation $\{ u_f^*,\hat{y}_f(u_f^*,w^*) \}$ can be far from the ground truth and results in a lack of robustness against unknown discrepancy with the SPC predictor.}

\subsection{Min-Max Robust DDPC with Guarantees}
The above analysis motivates us to develop a new robust DDPC formulation by adopting a \textit{pessimistic} min-max approach:
\textcolor{black}{
\begin{equation}
    \label{equation: robust DDPC with constraint}
    \begin{aligned}
        \min_{u_f\in\mathbb{U}}~&\max_{w\in\mathcal{W}(\Lambda)}\mathcal{J}(u_f,\hat{y}_f(u_f,w))\\
        {\rm s.t.}~&~~\hat{y}_f(u_f,w)\in\mathbb{Y},~\forall w\in\mathcal{W}(\Lambda).
    \end{aligned}
\end{equation}
This min-max formulation represents a fundamental paradigm shift from the min-min strategy used in PBR-DDPC \eqref{equation: regularized DDPC with constraint}, where $w$ is treated as an auxiliary decision variable to seek a best-case solution. In contrast, our formulation rigorously treats $w$ as a bounded certainty to optimize the worst-case scenario. The robustness of the proposed approach is twofold. In terms of control performance, the inner optimization actively identifies the adversarial realization of $w$ within the uncertainty set $\mathcal{W}(\Lambda)$ that maximizes the inner objective, which is subsequently minimized by the outer loop. In terms of constraint satisfaction, the output prediction is constrained to strictly reside within the feasible set $\mathbb{Y}$ for all possible realizations of $w$ inside the uncertainty set $\mathcal{W}(\Lambda)$, which prevents constraint violations under any realizations of $w$.} Consequently, the min-max formulation \eqref{equation: robust DDPC with constraint} can be seen as a ``robustified'' version of SPC, since the optimal control sequence $u_f^*$ derived is endowed with desirable robustness against the unknown difference between the output predictions from the SPC predictor and the realized output trajectory.

% As opposed to the min-min formulation \eqref{equation: regularized DDPC with constraint} in which $w$ is a decision variable, our min-max formulation \eqref{equation: robust DDPC with constraint} treats $w$ as uncertainty and thus bears a different philosophy. Amongst all output predictions $\hat{y}_f$ that result from $w\in\mathcal{W}(\Lambda)$ and deviate from the predictions from SPC predictor, our formulation \eqref{equation: robust DDPC with constraint} pursues the optimal \textit{worst-case} performance in terms of control cost minimization and constraint satisfaction. Specifically, the inner optimization problem identifies the adversarial realization $w^*$ within the uncertainty set $\mathcal{W}(\Lambda)$ that maximizes the inner objective, and the derived worst-case cost function is further minimized in the outer loop. Regarding the constraint satisfaction, the hard constraint of $\hat{y}_f$ demands that the output prediction strictly resides within the feasible set $\mathbb{Y}$ for all possible realizations of $w$ inside the set $\mathcal{W}$. By contrast, the PBR-DDPC in \eqref{equation: regularized DDPC with constraint} merely focuses on the \textit{best case}. Thus, the min-max formulation \eqref{equation: robust DDPC with constraint} can be seen as a ``robustified'' version of SPC, since the optimal control sequence $u_f^*$ derived is endowed with desirable robustness against the unknown difference between the output predictions from the SPC predictor and the realized output trajectory.

The conservatism of the min-max formulation \eqref{equation: robust DDPC with constraint} can be readily manipulated using a single parameter $\Lambda > 0$. We emphasize that $w$ in \eqref{equation: robust DDPC with constraint} arises from the non-uniqueness of solution $g$, and thus does not represent real uncertainty in practice, such as random noise. \textcolor{black}{Notably, in the noise-free case, the optimal input sequence and the optimal cost of the min-max formulation \eqref{equation: robust DDPC with constraint} coincide with those from SPC \eqref{equation: SPC} because the relation $Y_f\Phi^\bot=0$ strictly holds in the absence of noise, rendering the vector $w$ and the uncertainty set $\mathcal{W}(\Lambda)$ inactive.}

Next, we establish theoretical guarantees of our min-max formulation \eqref{equation: robust DDPC with constraint} in handling bounded noise added to input and output data. \textcolor{black}{For brevity, we omit the variables $u_f$ and $w$ in $\hat{y}_f(u_f,w)$ when no confusion arises.}

\begin{ass}\label{ass: 1}
All input and output measurements collected from \eqref{equation: innovation-form LTI system} are inexact and corrupted by bounded additive noise:
\begin{equation}
    u(t)=\bar{u}(t)+\xi_u(t),~ y(t)=\bar{y}(t)+\xi_y(t),
\end{equation}
where $\|\xi_u(t)\|\le\bar{\xi}_u$ and 
 $\|\xi_y(t)\|\le\bar{\xi}_y$.
\end{ass}

Indeed, Assumption \ref{ass: 1} is more realistic and has been widely adopted in DDPC literature \cite{berberich2021data,van2022noisy}. It brings uncertainty into $L$-long input-output trajectories in a behavioral setting, where an important concept for robust control design is data consistency. 
\begin{defi}[Consistent output data trajectories]
\label{definition: consistent output data trajectories}
The set of all plausible future output trajectories is defined as those consistent with given past data $\{u_p,y_p\}$ and future input $u_f$ for \textit{some} instance of $\xi_u(t)$ and $\xi_y(t)$, i.e.,
\begin{equation}
\label{equation: set of all plausible future output trajectories}
\begin{split}
    &\mathcal{Y}_o(u_p,y_p,u_f)\triangleq \\
    &\left\{ y_f \left| \begin{split}
 &{\rm col}(
                    \phi(u_f),y_f
                ) = {\rm col}(
                    \bar{\phi}(\bar{u}_f) + \xi_\phi,\bar{y}_f+\xi_{y_f}
                ),\\
                &{\rm col}(
                    \bar{\phi}(\bar{u}_f),\bar{y}_f
                )~ \text{is a trajectory of (1)}, \\ & \|\xi_u(t)\|\le\bar{\xi}_u,~\|\xi_{y}(t)\|\le\bar{\xi}_y,~ \forall t    
    \end{split} \right . \right \},            
\end{split}
\end{equation}
where $\xi_\phi\triangleq{\rm col}(\xi_{u_p},\xi_{u_f},\xi_{y_p})$, and $\xi_{u_p}$, $\xi_{u_f}$, $\xi_{y_p}$ and $\xi_{y_f}$ denote noise trajectories added to $\bar{u}_p$, $\bar{u}_f$, $\bar{y}_p$ and $\bar{y}_f$, respectively.
\end{defi}

By Assumption \ref{ass: 1}, Hankel data matrices in the min-max formulation \eqref{equation: robust DDPC with constraint} are also corrupted, denoted as $\Phi=\bar{\Phi}+\Xi_{\Phi},~Y_f=\bar{Y}_f+\Xi_Y$. To describe their impact on output prediction, we define the set of all possible output predictions subject to uncertainty $w$ as:
\begin{equation}
    \label{equation: definition of hatY}
    \begin{aligned}
        &\hat{\mathcal{Y}}_o(u_p,y_p,u_f)\\
        \triangleq&\left\{\hat{y}_f\left|\hat{y}_f=Y_f\Phi^\dagger\phi(u_f)+Y_f\Phi^\bot w,~w\in\mathcal{W}(\Lambda)\right.\right\}.
    \end{aligned}
\end{equation}
In fact, it is non-trivial to attain an explicit parameterization of $\mathcal{Y}_o(u_p,y_p,u_f)$ since exact data are practically inaccessible. Based on the link between $\mathcal{Y}_o(u_p,y_p,u_f)$ and $\hat{\mathcal{Y}}_o(u_p,y_p,u_f)$, the following result presents the expression of a sufficiently large but finite $\Lambda$ such that the set \eqref{equation: definition of hatY} covers actual noisy future outputs under Assumption \ref{ass: 1}.
\begin{theorem}
    \label{Theorem: bounded Lambda}
    Suppose that Assumption \ref{ass: 1} holds and $Y_f\Phi^\bot$ has full row rank. Given input-output data $\{u_p,y_p,u_f\}$, if $\Lambda$ is chosen as
    \begin{equation}
        \label{equation: condition of Lambda}
        \Lambda \ge \Lambda_o=\left\|\Phi^\bot(Y_f\Phi^\bot)^\dagger\right\|^2\left(\Lambda_1\left\|\phi(u_f)\right\|+\Lambda_2\right)^2,
    \end{equation}
    where
    \begin{equation}
        \label{equation: definition of Lambda12}
        \begin{aligned}
            \Lambda_1&=\sqrt{T}\left(\tfrac{1+\sqrt{5}}2\beta^2\xi_1\left\|Y_f\right\|+\delta\xi_2\right),\\
            \Lambda_2&=\delta\left(\left\|Y_f\right\|+\xi_2\sqrt{T}\right)\xi_1+\xi_2,
        \end{aligned}
    \end{equation}
    with $\delta = 1/ \sigma_{n_uL+n_x}(\bar{\Phi}),~\beta=\max\left\{\delta,~\left\|\Phi^\dagger\right\|\right\},~\xi_1=\bar{\xi}_u\sqrt{n_uL}+\bar{\xi}_y\sqrt{n_yL_p},~\xi_2=\bar{\xi}_y\sqrt{n_yL_f}$, then $\hat{\mathcal{Y}}_o(u_p,y_p,u_f)$ is a safe approximation of $\mathcal{Y}_o(u_p,y_p,u_f)$, i.e., $\mathcal{Y}_o(u_p,y_p,u_f)\subseteq\hat{\mathcal{Y}}_o(u_p,y_p,u_f)$.
\end{theorem}
\begin{proof}
    Consider $y_f\in\mathcal{Y}_o(u_p,y_p,u_f)$ and define
    \begin{equation}
        \label{equation: exists a w}
        \tilde{w}=(Y_f\Phi^\bot)^\dagger\underbrace{\left(y_f-Y_f\Phi^\dagger\phi(u_f)\right)}_{\triangleq\tilde{y}_f},
    \end{equation}
    which satisfies:
    \begin{equation}
        \label{equation: true output with noisy data}
        \begin{aligned}
            y_f&=Y_f\Phi^\dagger\phi(u_f)+Y_f\Phi^\bot\tilde{w},
        \end{aligned}
    \end{equation}
    due to the full row rank of $Y_f\Phi^\bot$. Thus, to prove $y_f\in\hat{\mathcal{Y}}_o(u_p,y_p,u_f)$, it suffices to show $\tilde{w}\in\mathcal{W}(\Lambda)$, which amounts to $\|\Phi^\bot\tilde{w}\|^2\le\Lambda$. Based on \eqref{equation: exists a w}, we have
    \begin{equation}
        \label{equation: decomposition of tildeyf}
        \|\Phi^\bot\tilde{w}\|=\|\Phi^\bot(Y_f\Phi^\bot)^\dagger\tilde{y}_f\|\le\|\Phi^\bot(Y_f\Phi^\bot)^\dagger\|\|\tilde{y}_f\|.
    \end{equation}
    According to Definition \ref{definition: consistent output data trajectories}, $y_f$ can be expressed as $y_f=\bar{y}_f+\xi_{y_f}=\bar{Y}_f\bar{\Phi}^\dagger\bar{\phi}(\bar{u}_f)+\xi_{y_f}$. It then follows that
    \begin{equation}
        \label{equation: upper bound decomposition}
        \begin{aligned}
             & \left\|\tilde{y}_f\right\|=\left\|\bar{Y}_f\bar{\Phi}^\dagger\bar{\phi}(\bar{u}_f)+\xi_{y_f}-Y_f\Phi^\dagger\phi(u_f)\right\|                                                          \\
              =&\left\|\left(\bar{Y}_f\bar{\Phi}^\dagger-Y_f\Phi^\dagger\right)\phi(u_f)-\bar{Y}_f\bar{\Phi}^\dagger\xi_\phi+\xi_{y_f}\right\|                                                   \\
              \le&\left\|\bar{Y}_f\bar{\Phi}^\dagger-Y_f\Phi^\dagger\right\|\left\|\phi(u_f)\right\|+\left\|\bar{Y}_f\bar{\Phi}^\dagger\right\|\left\|\xi_\phi\right\|+\left\|\xi_{y_f}\right\|.
        \end{aligned}
    \end{equation}    
    Firstly, we deal with the term $\left\|\bar{Y}_f\bar{\Phi}^\dagger-Y_f\Phi^\dagger\right\|$. Using the relation $Y_f=\bar{Y}_f+\Xi_Y$, we have
    \begin{equation}
        \label{equation: upper bound first part}
        \begin{aligned}
            \left\|\bar{Y}_f\bar{\Phi}^\dagger-Y_f\Phi^\dagger\right\| & =\left\|Y_f\left(\bar{\Phi}^\dagger-\Phi^\dagger\right)-\Xi_Y\bar{\Phi}^\dagger\right\|                                  \\
                                                                       & \le\left\|Y_f\right\|\left\|\bar{\Phi}^\dagger-\Phi^\dagger\right\|+\left\|\Xi_Y\right\|\left\|\bar{\Phi}^\dagger\right\|.
        \end{aligned}
    \end{equation}
    It follows from \cite[The general theorem]{stewart1977perturbation} that:
    \begin{equation}
        \label{equation: inverse bound}
        \begin{aligned}
            \left\|\bar{\Phi}^\dagger-\Phi^\dagger\right\| & \le\tfrac{1+\sqrt{5}}2\max\left\{\left\|\bar{\Phi}^\dagger\right\|^2,\left\|\Phi^\dagger\right\|^2\right\}\left\|\Xi_\Phi\right\| \\
                                                           & \le\tfrac{1+\sqrt{5}}2\beta^2\left\|\Xi_\Phi\right\|.
        \end{aligned}
    \end{equation}
    Since ${\rm rank}(\bar{\Phi})=n_uL+n_x$, $\sigma_{n_uL+n_x}(\bar{\Phi})$ is exactly the nonzero minimum singular value of $\bar{\Phi}$, from which the second inequality in \eqref{equation: inverse bound} follows. By definition of $\delta$, it holds that $\left\|\bar{\Phi}^\dagger\right\|\le\delta$. As a result, \eqref{equation: upper bound first part} can be further bounded by
    \begin{equation}
        \label{equation: upper bound first part complete}
        \begin{aligned}
            \left\|\bar{Y}_f\bar{\Phi}^\dagger-Y_f\Phi^\dagger\right\| & \le\tfrac{1+\sqrt{5}}2\beta^2\left\|Y_f\right\|\left\|\Xi_\Phi\right\|+\delta\left\|\Xi_Y\right\|.
        \end{aligned}
    \end{equation}

    Then we turn to the second term on the right-hand side of the inequality in \eqref{equation: upper bound decomposition}. For $\left\|\bar{Y}_f\bar{\Phi}^\dagger\right\|$, we have
    \begin{equation}
        \label{equation: upper bound second part}
        \left\|\bar{Y}_f\bar{\Phi}^\dagger\right\|=\left\|(Y_f+\Xi_Y)\bar{\Phi}^\dagger\right\|\le\left(\left\|Y_f\right\|+\left\|\Xi_Y\right\|\right)\delta.
    \end{equation}
    As for noise terms in data vectors and matrices, by utilizing the bounds $\{\bar{\xi}_u,\bar{\xi}_y\}$ of additive noise in Assumption \ref{ass: 1}, we have
    \begin{equation}
        \label{equation: upper bound of uncertainties}
        \begin{aligned}
             & \left\|\Xi_\Phi\right\|\le\bar{\xi}_u\sqrt{n_uLT}+\bar{\xi}_y\sqrt{n_yL_pT}=\xi_1\sqrt{T},                               \\
             & \left\|\Xi_Y\right\|\le\bar{\xi}_y\sqrt{n_yL_fT}=\xi_2\sqrt{T},~\left\|\xi_{y_f}\right\|\le\bar{\xi}_y\sqrt{n_yL_f}=\xi_2,           \\
             & \left\|\xi_\phi\right\|\le\bar{\xi}_u\sqrt{n_uL}+\bar{\xi}_y\sqrt{n_yL_p}=\xi_1.
        \end{aligned}
    \end{equation}
    By combining \eqref{equation: decomposition of tildeyf}, \eqref{equation: upper bound decomposition}, \eqref{equation: upper bound first part complete}, \eqref{equation: upper bound second part} and \eqref{equation: upper bound of uncertainties}, it follows that $\|\Phi^\bot\tilde{w}\|^2\le\Lambda_o\le\Lambda$, which indicates that  $y_f\in\hat{\mathcal{Y}}_o(u_p,y_p,u_f)$ and thus yields the desired result.
\end{proof}
{\color{black}
\begin{remark}
    It is not restrictive to assume the full row rank of $Y_f\Phi^\bot$ when the input-output data are corrupted by additive noise. Indeed, this is trivially satisfied in the presence of random noise, and can be evidenced in various empirical experiments. More importantly, because both $Y_f$ and $\Phi$ can be constructed from the known input-output data, the matrix $Y_f\Phi^\bot$ is accessible in practice, and its rank/singular values can be used to quantify the impact of noise and verify the rank condition in Theorem \ref{Theorem: bounded Lambda}.
\end{remark}
}

Theorem \ref{Theorem: bounded Lambda} clarifies when the case of additive bounded noise can be effectively encapsulated by the proposed min-max formulation using $\mathcal{W}(\Lambda)$ to describe uncertainty, despite some conservatism of \eqref{equation: condition of Lambda}. Note that Theorem \ref{Theorem: bounded Lambda} cannot be used as a guideline to choose $\Lambda$ for the uncertainty set. On the one hand, its theoretical bound is conservative; on the other hand, some values in \eqref{equation: condition of Lambda} are unknown, such as the parameter $\delta$ and the upper bounds $\{\bar{\xi}_u,\bar{\xi}_y\}$ of additive noise, which prevents the bound in Theorem \ref{Theorem: bounded Lambda} from being accessible in practice. In general, as $\{\bar{\xi}_u,\bar{\xi}_y\}$ increase, a larger $\Lambda$ is needed to describe the induced output prediction errors. In practice, $\Lambda$ serves as a flexible tuning parameter that can be determined by trial-and-error. In principle, $\Lambda$ should be sufficiently large such that the set $\hat{\mathcal{Y}}_o$ can encompass the noisy future outputs, while being as small as possible to avoid over-conservatism. This implies a simple and empirical strategy to select $\Lambda$, which solves the following optimization problem given a finite input-output sequence $\{u_p^i,u_f^i,y_p^i,y_f^i\}_{i=1}^{N_e}$:
\begin{equation}
    \label{equation: empirically calculate Lambda}
    \begin{aligned}
        \Lambda=\min_{\Lambda}&~\Lambda\\
        {\rm s.t.}&~y_f^i\in\hat{\mathcal{Y}}_o(u_p^i,u_f^i,y_p^i),~i\in[N_e],
    \end{aligned}
\end{equation}
where $\hat{\mathcal{Y}}_o$ defined in \eqref{equation: definition of hatY} is implicitly dependent on $\Lambda$.

% It is seen that $\Lambda_o$ in \eqref{equation: condition of Lambda} only relies on the collected noisy input-output data $\{u^d(t),y^d(t)\}_{t=1}^N$, the online past data $\{u_p,y_p\}$, the parameters $\{L_p,L_f,N\}$ of the DDPC method, and the upper bounds $\{\bar{\xi}_u,\bar{\xi}_y\}$ of additive noise. As $\{\bar{\xi}_u,\bar{\xi}_y\}$ increase, a larger $\Lambda$ is needed to describe the induced output prediction errors. 

% When there is no noise in input-output data, i.e., $\bar{\xi}_u=\bar{\xi}_y=0$, $\Lambda$ can be set to zero according to \eqref{equation: condition of Lambda} and thus \eqref{equation: robust DDPC with constraint} reduces to the problem of SPC \eqref{equation: SPC}.

% Theorem \ref{theorem: relation between realized cost and optimal cost} shows that solving the min-max problem \eqref{equation: robust DDPC with constraint} yields a robustified control sequence that intrinsically encapsulates the effect of bounded additive noise. 

\textcolor{black}{Furthermore, since the actual future output $y_f$ is guaranteed to reside within the safe approximation set $\hat{\mathcal{Y}}_o$ under certain conditions of $\Lambda$, one can readily establish through the triangle inequality that the true control cost $\mathcal{J}(u_f^*,y_f)$ is deterministically upper-bounded by a term depending on $\mathcal{J}(u^*_f,\hat{y}_f^*)$ and $\Lambda$ \cite[Theorem 3.1]{huang2023robust}.}

\subsection{Tractable Reformulation for Robust DDPC}
For DDPC algorithms, it is crucial to ensure the tractability of optimization problems and the efficiency of online solutions. \textcolor{black}{Note that the uncertainty vector $w\in\mathbb{R}^T$ resides in a high-dimensional space that grows with sample size $N$. Nevertheless, it affects the min-max problem solely through $Y_f\Phi^\bot w$. To eliminate redundant degrees of freedom without introducing conservatism, we perform a compact singular value decomposition (SVD) as $Y_f\Phi^\bot=W\Sigma V^\top$, where $\Sigma\in\mathbb{R}^{n_z\times n_z}$ and $n_z={\rm rank}(Y_f\Phi^\bot)\le n_yL_f$. By defining a reduced-dimensional uncertainty variable $z=V^\top w\in\mathbb{R}^{n_z}$ and a matrix $M_z=W\Sigma$, the prediction error term $Y_f\Phi^\bot w$ becomes exactly $M_zz$. Since $V$ spans a subspace of ${\rm im}(\Phi^\bot)$, we have $\Phi^\bot V=V$. It follows that $\|z\|_2^2=\|V^\top\Phi^\bot w\|_2^2\le\|\Phi^\bot w\|_2^2\le\Lambda$. Conversely, for any $z$ satisfying $\|z\|^2_2\le\Lambda$, there always exists a valid realization $w=Vz\in\mathcal{W}(\Lambda)$ that achieves it. Consequently, the uncertainty set $\mathcal{W}(\Lambda)$ can be strictly and equivalently redefined in a much lower-dimensional space as $\mathcal{Z}(\Lambda)=\left\{z\in\mathbb{R}^{n_z}~|~\|z\|^2\le\Lambda\right\}$. In this context, the original min-max problem \eqref{equation: robust DDPC with constraint} can be reformulated into a more computationally efficient form as:
\begin{equation}
    \label{equation: robust DDPC with constraint with z}
    \begin{aligned}
        \min_{u_f\in\mathbb{U}}&~\max_{z\in\mathcal{Z}(\Lambda)}~\mathcal{J}(u_f,\hat{y}_f(u_f,z))\\
        {\rm s.t.}&~~\hat{y}_f(u_f,z)\in\mathbb{Y},~\forall z\in\mathcal{Z}(\Lambda),
    \end{aligned}
\end{equation}
where the function of output prediction is defined as $\hat{y}_f(u_f,z)=Y_f\Phi^\dagger\phi(u_f)+M_zz$.}

In what follows, we show how to convert the min-max problem \eqref{equation: robust DDPC with constraint with z} into a convex program that is amenable to off-the-shelf numerical solvers.
\begin{proposition}
    \label{proposition: robust DeePC}
    The input sequence $u_f^*$ is a minimizer of \eqref{equation: robust DDPC with constraint with z} if and only if there exists $(b^*,\psi^*,\gamma^*,\{\mu_i^*\}_{i=1}^{l_y})$ such that $u_f^*$ also minimizes the following SDP problem:
    \begin{subequations}
        \label{equation: robust DDPC semidefinite}
        \begin{align}
             \min_{\begin{subarray}{c}u_f,b,\psi,\\\gamma,\{\mu_i\}_{i=1}^{l_y}\end{subarray}}~&\psi+\|u_f\|_R^2   \label{equation: robust DDPC semidefinite: obj}               \\
             {\rm s.t.}\quad &
             \label{equation: robust DDPC semidefinite: b}
             u_f\in\mathbb{U},~b=Y_f\Phi^\dagger\phi(u_f),\\
            \label{equation: robust DDPC semidefinite: yf}
             & \begin{bmatrix}
                   1-\mu_i\Lambda & \star             & \star \\
                   0          & \mu_i I & \star       \\
                   G_y^i b+c_y^i         & G_y^iM_z          & I
               \end{bmatrix}\succeq0,~i\in[l_y],                                                               \\
            \label{equation: robust DDPC semidefinite: LMI}
             & \begin{bmatrix}
                   \psi-\gamma\Lambda+y_r^\top Q(2b-y_r) & \star     & \star \\
                   M_z^\top Q y_r                                      & \gamma I & \star \\
                   b                                       & M_z           & Q^{-1}
               \end{bmatrix}\succeq0,\\
               &\gamma\ge0,~\mu_i\ge0,~i\in[l_y].
        \end{align}
    \end{subequations}
    Moreover, the minima of \eqref{equation: robust DDPC with constraint with z} and \eqref{equation: robust DDPC semidefinite} coincide.
\end{proposition}
\begin{proof}
    \textcolor{black}{
    The proof is obtained by utilizing the S-Lemma and Schur complement, which are widely used in robust optimization \cite{elghaoui1997robust} and robust model predictive control (MPC) \cite{bekiroglu2013robust}. By defining $b$ in \eqref{equation: robust DDPC semidefinite}, the data-driven predictor in \eqref{equation: robust DDPC with constraint with z} can be succinctly expressed as $\hat{y}_f(u_f,z)=M_zz+b$. By introducing an auxiliary variable $\psi$, the min-max problem \eqref{equation: robust DDPC with constraint with z} can be rewritten as:
    \begin{subequations}
        \label{equation: robust DDPC with constraint and t}
        \begin{align}
            \min_{u_f,\psi,b,\{\mu_i\}_{i=1}^{l_y}} & ~\psi+\|u_f\|_R^2\\
            {\rm s.t.}~~\quad                       & u_f\in\mathbb{U},~b=Y_f\Phi^\dagger\phi(u_f),\\
            \label{equation: robust DDPC with constraint and t: output}
            &\hat{y}_f(u_f,z)\in\mathbb{Y},~\forall z\in\mathcal{Z}(\Lambda), \\
            \label{equation: robust DDPC with constraint and t: max}
            & \psi\ge\|M_zz+b-y_r\|_Q^2,~\forall z\in\mathcal{Z}(\Lambda),
        \end{align}
    \end{subequations}
    where the inner maximization problem in \eqref{equation: robust DDPC with constraint with z} is equivalently converted into the infinite constraint \eqref{equation: robust DDPC with constraint and t: max}. By invoking the S-Lemma \cite[Theorem B.2.1]{ben2009robust} to individually tackle the output constraints in \eqref{equation: robust DDPC with constraint and t: output}, the constraints hold if and only if there exist multipliers $\{\mu_i\ge0\}_{i=1}^{l_y}$ such that
    \begin{equation}
        \label{equation: robust DDPC with constraint: yf inequality}
        \begin{aligned}
            \begin{bmatrix}
                -\mu_i\Lambda-\|G_y^ib+c_y^i\|^2+1 & -(G_y^ib+c_y^i)^\top G_y^iM_z\\
                -(G_y^iM_z)^\top(G_y^ib+c_y^i) & \mu_iI-(G_y^iM_z)^\top(G_y^iM_z)
            \end{bmatrix}\succeq0,&\\
            i\in[l_y].&
        \end{aligned}
    \end{equation}
    Using the Schur complement \cite{cottle1974manifestations}, the above LMIs then become \eqref{equation: robust DDPC semidefinite: yf}.
    The constraint \eqref{equation: robust DDPC with constraint and t: max} can be tackled in a similar spirit to \eqref{equation: robust DDPC with constraint and t: output}. Invoking the S-Lemma, the constraint \eqref{equation: robust DDPC with constraint and t: max} holds if and only if there exists $\gamma\ge0$ such that
    \begin{equation*}
        % \label{equation: equation: robust DDPC with constraint and t: max: matrix form}
        \begin{aligned}
            \begin{bmatrix}
                -\gamma\Lambda-(b-y_r)^\top Q(b-y_r)+\psi & -(b-y_r)^\top QM_z\\
                -M_z^\top Q(b-y_r) & \gamma I-M_z^\top QM_z
            \end{bmatrix}\succeq0,
        \end{aligned}
    \end{equation*}
    thereby giving rise to \eqref{equation: robust DDPC semidefinite: LMI} by the Schur complement. Up to now, we have arrived at the tractable reformulation \eqref{equation: robust DDPC semidefinite}, which completes the proof.}
    % One can further rewrite \eqref{equation: robust DDPC with constraint: yf inequality} as the following LMIs
    % \begin{equation}
    %     \label{equation: robust DDPC with constraint: yf S-Lemma}
    %     \begin{aligned}
    %           &\begin{bmatrix}
    %                1-\mu_i\Lambda & 0 \\0&\mu_i I
    %            \end{bmatrix}\\    
    %           &-\begin{bmatrix}
    %                 (G_y^ib+c_y^i)^\top \\(G_y^iM_z)^\top
    %             \end{bmatrix} \begin{bmatrix}
    %                                    G_y^ib+c_y^i & G_y^iM_z
    %                                \end{bmatrix}\succeq0,~i\in[l_y].
    %     \end{aligned}
    % \end{equation}

    % which amounts to
    % \begin{equation}
    %     \begin{aligned}
    %         & \begin{bmatrix}
    %                \psi-\gamma\Lambda+y_r^\top Q(2b-y_r) & y_r^\top QM_z \\ (y_r^\top QM_z)^\top& \gamma I
    %           \end{bmatrix} \\
    %         & \qquad \qquad \qquad  -\begin{bmatrix}
    %                 b^\top \\ M_z^\top
    %             \end{bmatrix}Q\begin{bmatrix}
    %                                     b & M_z
    %                                 \end{bmatrix} \succeq 0,
    %     \end{aligned}
    % \end{equation}
\end{proof}

\section{Feedback Robust DDPC}

In the previous section, an open-loop input sequence is produced by the proposed robust DDPC at every time step. However, this tends to be conservative, often leading to a severe underestimation of the set of feasible trajectories \cite{kerrigan2004feedback}. Inspired by feedback robust MPC \cite{lofberg2003approximations,li2010feedback}, we propose a new feedback robust DDPC scheme by incorporating an affine feedback control policy such that the control action currently applied is endowed with a better foresight into the future horizon. Specifically, we adopt the following feedback control policy, which is an affine function of output prediction error \cite{lofberg2003approximations,bemporad1998reducing}:
\begin{equation}
    \label{equation: affine decision rule with prediction error}
    u_f=v_f+K\left(y_f-\hat{y}_{f,v}\right),
\end{equation}
where $v_f\in\mathbb{R}^{n_uL_f}$ is a nominal input, $K\in\mathbb{R}^{n_uL_f\times n_yL_f}$ is a feedback gain, and $\hat{y}_{f,v}=Y_f\Phi^\dagger\phi(v_f)$ is the ``nominal'' output prediction by SPC predictor with $v_f$ used as the future input trajectory. To ensure causality of the affine feedback policy in \eqref{equation: affine decision rule with prediction error}, the feedback gain matrix $K$ should have a strictly lower-block triangular structure, where all entries along and above the diagonal $n_y\times n_u$ blocks are enforced to be zero. 

In the min-max formulation \eqref{equation: robust DDPC with constraint}, the prediction error of SPC is assumed to be $Y_f\Phi^\bot w$. Following the spirit in \eqref{equation: affine decision rule with prediction error}, we replace the input trajectory $u_f$ in \eqref{equation: robust DDPC with constraint} with the following affine feedback control policy:
\textcolor{black}{
\begin{equation}
    \label{equation: affine decision rule}
    u_f(v_f,K,w)=v_f+KY_f\Phi^\bot w,
\end{equation}
}
where $Y_f\Phi^\bot w$ is designed as the feedback term corresponding to the output prediction error, \textcolor{black}{and the future input is defined as a function of $\{v_f,K,w\}$ to emphasize its dependence on these variables.} In this case, the output prediction under the affine feedback policy \eqref{equation: affine decision rule} becomes:
\textcolor{black}{
\begin{equation}
    \label{equation: OP with feedback control law}
    \begin{aligned}
        \hat{y}_f(v_f,K,w) & =Y_f\Phi^\dagger{\rm col}(
                                        u_p,KY_f\Phi^\bot w+v_f,y_p
                                    )+Y_f\Phi^\bot w                    \\
                  & =Y_f\Phi^\dagger\phi(v_f)+\left(I+M_fK\right)Y_f\Phi^\bot w,
    \end{aligned}
\end{equation}
}
where $M_f\in\mathbb{R}^{n_yL_f \times n_uL_f}$ is derived from
\begin{equation}
    \label{equation: M1 and M2}
    Y_f\Phi^\dagger\phi(u_f)=M_fu_f+M_p{\rm col}(u_p,y_p).
\end{equation}
Based on \eqref{equation: affine decision rule}, \eqref{equation: OP with feedback control law}, and \eqref{equation: robust DDPC with constraint}, a feedback robust version of DDPC can be derived as the following min-max optimization problem:
\textcolor{black}{
\begin{equation}
    \label{equation: robust DDPC with optimized feedback gain}
    \begin{aligned}
        \min_{v_f,K}~&\max_{w\in\mathcal{W}(\Lambda)}~\mathcal{J}\left(u_f(v_f,K,w),\hat{y}_f(v_f,K,w)\right)\\
        {\rm s.t.}~&~~u_f(v_f,K,w)\in\mathbb{U},~\forall w\in\mathcal{W}(\Lambda),\\
        &~~\hat{y}_f(v_f,K,w)\in\mathbb{Y},~\forall w\in\mathcal{W}(\Lambda),\\
        &~~K~\textrm{strictly lower-block triangular},
    \end{aligned}
\end{equation}
}
where the nominal input $v_f$ and feedback gain $K$ are jointly optimized as decision variables. As compared to the open-loop min-max formulation \eqref{equation: robust DDPC with constraint}, the feedback version \eqref{equation: robust DDPC with optimized feedback gain} allows future controls to be adapted to uncertainty in dynamical systems. It is noted that both input and output constraints shall be robustly satisfied in \eqref{equation: robust DDPC with optimized feedback gain} for all possible realizations of $w$. Besides, even if an affined feedback policy is produced by solving \eqref{equation: robust DDPC with optimized feedback gain}, a receding horizon implementation shall be carried out. That is, after solving for $v_f^*$ and $K^*$, only the first $n_u$ entries in $v_f^*$ shall be applied to the system as the input, and this procedure will be repeated when moving on to the next time instance.

Next, we establish the performance guarantee of our feedback robust DDPC scheme. The following definition is made, which is a generalization of Definition \ref{definition: consistent output data trajectories} to the case of feedback control policy \eqref{equation: affine decision rule with prediction error}.
\begin{defi}[Consistent output data trajectories in closed-loop]
    \label{definition: consistent output data trajectories under feedback policy}
    Given past data $\{u_p,y_p\}$ and parameters $\{v_f,K\}$ of feedback policy, the set of all plausible future output trajectories that are consistent with the given data for some instances of $\xi_u(t)$ and $\xi_y(t)$, i.e.,
    \begin{equation}
    \label{equation: set of all plausible future output trajectories under feedback}
    \begin{split}
        &\mathcal{Y}_c(u_p,y_p,v_f,K)\triangleq \\
        &\left\{ y_f \left| \begin{split}
     &{\rm col}(\phi(u_f),y_f) = {\rm col}(\bar{\phi}(\bar{u}_f)+\xi_\phi,\bar{y}_f+\xi_{y_f}),\\
     &u_f=v_f+K\left(y_f-\hat{y}_{f,v}\right),~\hat{y}_{f,v}=Y_f\Phi^\dagger\phi(v_f),\\
     &{\rm col}(\bar{\phi}(\bar{u}_f),\bar{y}_f)~\text{is a trajectory of (1)}, \\ 
                    & \|\xi_u(t)\|\le\bar{\xi}_u,~\|\xi_{y}(t)\|\le\bar{\xi}_y    
        \end{split} \right . \right \}.             
    \end{split}
    \end{equation}
\end{defi}
Here in Definition \ref{definition: consistent output data trajectories under feedback policy}, the affine feedback policy \eqref{equation: affine decision rule with prediction error} is considered instead of a deterministic input sequence. Meanwhile, we define the set of all output predictions that are considered in \eqref{equation: robust DDPC with optimized feedback gain} as:
\begin{equation}
    \label{equation: set of all plausible output predictions under feedback}
    \begin{aligned}
        &\hat{\mathcal{Y}}_c(u_p,y_p,v_f,K)\\
        \triangleq&\left\{\hat{y}_f\left|\hat{y}_f=Y_f\Phi^\dagger\phi(v_f)+(I+M_fK)Y_f\Phi^\bot w,~w\in\mathcal{W}(\Lambda)\right.\right\}.
    \end{aligned}
\end{equation}
By establishing a link between $\mathcal{Y}_c(u_p,y_p,v_f,K)$ and $\hat{\mathcal{Y}}_c(u_p,y_p,v_f,K)$, the following result presents the expression of a sufficiently large but finite $\Lambda$ such that the set \eqref{equation: set of all plausible future output trajectories under feedback} encompasses actual noisy future outputs in the noisy case.
\begin{theorem}
    \label{theorem: bound Lambda with feedback control policy}
    Let Assumption \ref{ass: 1} hold, and suppose that $Y_f\Phi^\bot$ has full row rank and $(I-\bar{M}_fK)(I+M_fK)$ is non-singular. Given past data $\{u_p,y_p\}$, nominal input $v_f$ and feedback gain $K$, if $\Lambda$ is chosen as
    \begin{equation}
        \label{equation: sufficient condition under feedback law}
        \begin{aligned}
            \Lambda\ge\Lambda_c=&\frac{\|\Phi^\bot(Y_f\Phi^\bot)\|^2}{\left\|(I-\bar{M}_fK)(I+M_fK)\right\|^2}
            \left(\Lambda_1\left\|\phi(v_f)\right\|+\Lambda_2\right)^2,
        \end{aligned}
    \end{equation}
    where $\Lambda_1$ and $\Lambda_2$ are defined as \eqref{equation: definition of Lambda12} and $\bar{M}_f$ is derived from
    \begin{equation}
        \bar{Y}_f\bar{\Phi}^\dagger\bar{\phi}(\bar{u}_f)=\bar{M}_f\bar{u}_f+\bar{M}_p{\rm col}(\bar{u}_p,\bar{y}_p),
    \end{equation}
    then $\hat{\mathcal{Y}}_c(u_p,y_p,v_f,K)$ is a safe approximation of $\mathcal{Y}_c(u_p,y_p,v_f,K)$, i.e., $\mathcal{Y}_c(u_p,y_p,v_f,K)\subseteq \hat{\mathcal{Y}}_c(u_p,y_p,v_f,K)$.
\end{theorem}
\begin{proof}
    Consider a future output sequence $y_f\in\mathcal{Y}_c(u_p,y_p,v_f,K)$ and define
    \begin{equation}
        \label{equation: exists a w with feedback policy}
        \begin{aligned}
            \tilde{w}=(Y_f\Phi^\bot)^{\dagger}(I+M_fK)^{-1}\underbrace{\left(y_f-Y_f\Phi^\dagger\phi(v_f)\right)}_{\triangleq\tilde{y}_{f,v}},
        \end{aligned}
    \end{equation}
    which satisfies $y_f=Y_f\Phi^\dagger\phi(v_f)+\left(I+M_fK\right)Y_f\Phi^\bot\tilde{w}$ owing to the full row rank of $Y_f\Phi^\bot$ and the non-singularity of $(I+M_fK)$. Thus, to prove $y_f\in\hat{\mathcal{Y}}_c(u_p,y_p,v_f,K)$, it suffices to show $\tilde{w}\in\mathcal{W}(\Lambda)$, which amounts to $\|\Phi^\bot\tilde{w}\|^2\le\Lambda_c$. By Definition \ref{definition: consistent output data trajectories under feedback policy}, $y_f$ can be expressed as:
    \begin{equation}
        \label{equation: true output with noise-free data under feedback law: pre}
        \begin{aligned}
            y_f&=\bar{y}_f+\xi_{y_f}=\bar{Y}_f\bar{\Phi}^\dagger{\rm col}(\bar{u}_p,u_f-\xi_{u_f},\bar{y}_p)+\xi_{y_f}.
        \end{aligned}
    \end{equation}
    Then, by utilizing \eqref{equation: affine decision rule with prediction error} and the invertibility of $(I-\bar{M}_fK)$, we have:
    \begin{equation}
        \label{equation: true output with noise-free data under feedback law}
        \begin{aligned}
             &y_f=\bar{Y}_f\bar{\Phi}^\dagger{\rm col}(                \bar{u}_p,v_f+K(y_f-\hat{y}_{f,v})-\xi_{u_f},\bar{y}_p)+\xi_{y_f}\\
             &y_f=\bar{Y}_f\bar{\Phi}^\dagger\bar{\phi}(\bar{v}_f)+\bar{M}_fK\left(y_f-Y_f\Phi^\dagger\phi(v_f)\right)+\xi_{y_f} \\
            &(I-\bar{M}_fK)y_f=\left(\bar{Y}_f\bar{\Phi}^\dagger\bar{\phi}(\bar{v}_f)-\bar{M}_fKY_f\Phi^\dagger\phi(v_f)+\xi_{y_f}\right),\\
            & y_f=(I-\bar{M}_fK)^{-1}\left(\bar{Y}_f\bar{\Phi}^\dagger
            \bar{\phi}(\bar{v}_f)-\bar{M}_fKY_f\Phi^\dagger
            \phi(v_f)+\xi_{y_f}\right),        
        \end{aligned}
    \end{equation}
    where $\bar{v}_f=v_f-\xi_{u_f}$ is similar to the relation between $u_f$ and $\bar{u}_f$ in Definition \ref{definition: consistent output data trajectories}. It then follows for $\tilde{y}_{f,v}$ that
    \begin{equation}
        \begin{aligned}
            \tilde{y}_{f,v}=&(I-\bar{M}_fK)^{-1}\left(\bar{Y}_f\bar{\Phi}^\dagger\bar{\phi}(\bar{v}_f)+\xi_{y_f}\right)\\
            &-\left[(I-\bar{M}_fK)^{-1}\bar{M}_fK+I\right]Y_f\Phi^\dagger\phi(v_f)\\
            =&(I-\bar{M}_fK)^{-1}\underbrace{\left(\bar{Y}_f\bar{\Phi}^\dagger\bar{\phi}(\bar{v}_f)-Y_f\Phi^\dagger\phi(v_f)+\xi_{y_f}\right)}_{\triangleq\tilde{y}_{v}}.
        \end{aligned}
    \end{equation}
    Thus, based on \eqref{equation: exists a w with feedback policy}, we have
    \begin{equation}
        \label{equation: decomposition feedback}
        \begin{aligned}
            \|\Phi^\bot\tilde{w}\|&=\|\Phi^\bot(Y_f\Phi^\bot)^\dagger[(I-\bar{M}_fK)(I+M_fK)]^{-1}\tilde{y}_v\|\\
            &\le\|[(I-\bar{M}_fK)(I+M_fK)]^{-1}\|\|\Phi^\bot(Y_f\Phi^\bot)^\dagger\|\|\tilde{y}_v\|.
        \end{aligned}
    \end{equation}
    Following the proof of Theorem \ref{Theorem: bounded Lambda}, we have:
    \begin{equation}
        \label{equation: tildeyv feedback}
        \|\tilde{y}_v\|\le\Lambda_1\phi(v_f)+\Lambda_2.
    \end{equation}
    Plugging \eqref{equation: tildeyv feedback} into \eqref{equation: decomposition feedback}, we obtain $\|\Phi^\bot\tilde{w}\|^2\le\Lambda_c\le\Lambda$, which implies that $y_f\in\hat{\mathcal{Y}}_c(u_p,y_p,v_f,K)$ and thus completes the proof of $\mathcal{Y}_c(u_p,y_p,v_f,K)\subseteq\hat{\mathcal{Y}}_c(u_p,y_p,v_f,K)$.
\end{proof}

As shown in Theorem \ref{theorem: bound Lambda with feedback control policy}, there exists a sufficiently large $\Lambda$ as \eqref{equation: sufficient condition under feedback law} such that the prediction error induced by bounded additive noise under the feedback policy \eqref{equation: affine decision rule with prediction error} can be fully encapsulated by the term $Y_f\Phi^\bot w$ with the uncertainty set $\mathcal{W}(\Lambda)$. Compared with Theorem \ref{Theorem: bounded Lambda}, the requirement for $\Lambda$ in Theorem \ref{theorem: bound Lambda with feedback control policy} may be relaxed if $\|(I-\bar{M}_fK)(I+M_fK)\|\ge1$. 
\textcolor{black}{
Building upon the SVD-based dimensionality reduction introduced in \eqref{equation: robust DDPC with constraint with z}, the output prediction error $Y_f\Phi^\bot w$ in the affine feedback policy can be equivalently replaced by $M_zz$, where $M_z=W\Sigma$ and $z\in\mathcal{Z}(\Lambda)$. Consequently, the feedback robust DDPC problem \eqref{equation: robust DDPC with optimized feedback gain} can be reformulated into the following more compact form:
\begin{equation}
    \label{equation: robust DDPC with optimized feedback gain with z}
    \begin{aligned}
        \min_{v_f,K}~&\max_{z\in\mathcal{Z}(\Lambda)}~\mathcal{J}\left(u_f(v_f,K,z),\hat{y}_f(v_f,K,z)\right)\\
        {\rm s.t.}~&~~u_f(v_f,K,z)\in\mathbb{U},~\forall z\in\mathcal{Z}(\Lambda),\\
        &~~\hat{y}_f(v_f,K,z)\in\mathbb{Y},~\forall z\in\mathcal{Z}(\Lambda),\\
        &~~K~\textrm{strictly lower-block triangular},
    \end{aligned}
\end{equation}
with $u_f(v_f,K,z)=v_f+KM_zz$ and $\hat{y}_f(v_f,K,z)=Y_f\Phi^\dagger\phi(v_f)+(I+M_fK)M_zz$. To solve this min-max problem, we can derive a tractable convex SDP problem following the same reformulation procedure in Proposition \ref{proposition: robust DeePC}. 
}

\section{Simulation Case Studies}

To demonstrate the effectiveness of the proposed methods, we consider a benchmark two-mass-spring-damper system, where the state-space representation is adopted from \cite{o’dwyer2023data} as
\begin{equation}
    \label{equation: system of two mass-spring-damper}
    \left \{ \begin{aligned}
        x(t+1) & =Ax(t)+B_u\left[u(t)+v_1(t)\right] \\
        y(t)   & =Cx(t)+B_vv_2(t),
    \end{aligned} \right .
\end{equation}
with system matrices
\begin{equation*}
    \label{equation: matrices of two mass-spring-damper}
    \begin{aligned}
         & A=\begin{bmatrix}
                 1                            & 0                        & \Delta_t                      & 0                         \\
                 0                            & 1                        & 0                             & \Delta_t                  \\
                 -\frac{k_1}{m_1}\Delta_t & \frac{k_1}{m_1}\Delta_t  & 1-\frac{b_1}{m_1}\Delta_t & \frac{b_1}{m_1}\Delta_t   \\
                 \frac{k_1}{m_2}\Delta_t      & -\frac{k_1+k_2}{m_2}\Delta_t & \frac{b_1}{m_2}\Delta_t       & 1-\frac{b_1+b_2}{m_2}\Delta_t
             \end{bmatrix}, \\
         & B_u=\begin{bmatrix}
                   0 & 0 & \frac{1}{m_1}\Delta_t & 0
               \end{bmatrix}^\top,~B_v = \begin{bmatrix}
                                             0.5 & 1 & 0.4 & 0.3
                                         \end{bmatrix}^\top,\\
        & C= I_4.
    \end{aligned}
\end{equation*}
The system parameters are chosen as $k_1=4,~k_2=4,~c_1=1.5,~c_2=2,~m_1=1.2,~m_2=2$, and the sampling time is $\Delta_t=0.1\si{s}$. \textcolor{black}{Input disturbance $v_1(t)$ and measurement noise $v_2(t)$ are modeled as auto-regressive processes $v_i(t)=0.5v_i(t-1)+e_i(t)$, where $e_1(t)\sim\mathcal{N}(0,0.01^2),e_2(t)\sim\mathcal{N}(0,0.19^2)$ (truncated to $\pm3\sigma$). The primary task is to track a square wave reference $y_r(t)$ with $N_{\rm test}=100$ for the position of the first mass.}

For offline data collection, a square wave with a period of $600$ time-steps and amplitude of $1$, contaminated by a zero-mean Gaussian noise sequence with variance $0.01$, is used as the persistently exciting input to collect an informative input-output data trajectory with length $N=600$. The control inputs are subject to the constraint $[-5,5]$, while the velocities of the two masses are restricted to $[-1.4,1.4]$. The control parameters of DDPC are set to $L_p=L_f=5$, $Q=1$, and $R=0.01$. We compare four receding horizon strategies:
\begin{enumerate}
    \item[(1)] \textbf{SPC}: The classical SPC scheme \eqref{equation: SPC} \cite{favoreel1999spc}.
    \item[(2)] \textbf{PBR-DDPC}: The PBR-DDPC method \cite{dorfler2023bridging} reformulated as the min-min problem \eqref{equation: regularized DDPC with constraint}.
    \item[(3)] \textbf{R-DDPC}: The proposed robust DDPC \eqref{equation: robust DDPC semidefinite}.
    \item[(4)] \textbf{FR-DDPC}: The proposed feedback robust DDPC derived from \eqref{equation: robust DDPC with optimized feedback gain with z} with affine feedback control policy.
\end{enumerate}
% All strategies are implemented in a one-step receding horizon fashion, where the resultant optimization problems are solved at each time instance and only the first optimized input is applied to the system. All computations are performed in MATLAB 2024a on a desktop computer with Inter(R) Core(TM) i9-14900K CPU  @ 3.2 GHz and 64GB RAM. 
All optimization problems are solved using Mosek 10.2 \cite{mosek}. To evaluate the control performance of the DDPC methods, the following index is used:
\begin{equation}
    \mathcal{J}_{\rm total}=\sum_{k=1}^{N_{\rm test}}||y(k) - y_r(k)||_Q^2  + \sum_{k=1}^{N_{\rm test}}||u(k)||_R^2,
\end{equation}
where $u(k)$ is the control input generated from each control strategy, and $y(t)$ is the corresponding noisy output from the realized system.

We run $100$ Monte Carlo simulations to evaluate the control performance. For PBR-DDPC, R-DDPC, and FR-DDPC, the hyperparameter $\Lambda$ is selected via grid search across $[10^{-10},10^6]$. \textcolor{black}{Table \ref{tab:performance} summarizes the mean and standard deviation of the output tracking error cost and the control effort cost for all methods. It can be seen that PBR-DDPC, whose optimal performance is achieved with a very small $\Lambda$, closely resembles SPC with large tracking errors and high variance. This confirms that the optimism inherent to PBR-DDPC is not beneficial for improving the performance under uncertainty. In contrast, the two proposed robust DDPC methods have the upper hand over SPC and PBR-DDPC with lower tracking errors. By incorporating the affine feedback policy, FR-DDPC successfully mitigates the conservatism of the open-loop design. It not only achieves the most accurate tracking but also stabilizes the control actions due to a smaller variance of control effort.}

\begin{table}[h]
\centering
\caption{\textcolor{black}{Detailed Performance Metrics over 100 Monte Carlo Simulations}}
\label{tab:performance}
\renewcommand{\arraystretch}{1.3}
\resizebox{\columnwidth}{!}{% 自动缩放以适应单栏宽度
\textcolor{black}{
\begin{tabular}{lcc}
\hline\hline
\textbf{Method} & \textbf{Tracking Error} ($\sum \|y-y_r\|_Q^2$) & \textbf{Control Effort} ($\sum \|u\|_R^2$) \\ \hline
SPC         & $3.50 \pm 0.55$ & $0.56 \pm 0.14$ \\
PBR-DDPC   & $3.50 \pm 0.55$ & $0.56 \pm 0.14$ \\
R-DDPC   & $2.67 \pm 0.27$ & $0.90 \pm 1.71$ \\
FR-DDPC  & $2.52 \pm 0.24$ & $1.21 \pm 0.31$ \\ \hline\hline
\end{tabular}%
}}
\end{table}

\begin{figure}[ht]
    \centering
    \includegraphics[width=0.4\textwidth]{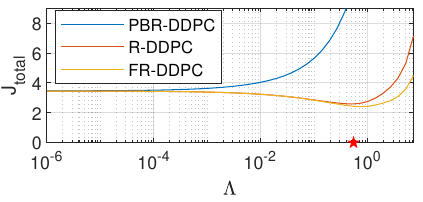}
    \caption{Average $\mathcal{J}_{\rm total}$ of different control algorithms under different $\Lambda$, where the pentagram denotes $\Lambda=0.5446$ obtained from \eqref{equation: empirically calculate Lambda}.}
    \label{fig: SNR30 Lambda}
\end{figure}

Then, we investigate the effect of the size parameter $\Lambda$ on the performance of PBR-DDPC, R-DDPC, and \textcolor{black}{FR-DDPC}. Fig. \ref{fig: SNR30 Lambda} depicts the average control costs in $100$ Monte Carlo simulations under varying values of $\Lambda$, where all three methods boil down to SPC when $\Lambda \to 0$. However, as $\Lambda$ grows, the performance of PBR-DDPC degrades severely due to its optimistic nature. Conversely, with $\Lambda$ suitably chosen (e.g. $10^{-2}<\Lambda<1$), the proposed min-max robust DDPC methods effectively improve upon SPC. Meanwhile, the value of $\Lambda$ obtained by the data-based rule \eqref{equation: empirically calculate Lambda} appears to be near-optimal, which justifies its use for determining $\Lambda$. 

% gets improved with $\Lambda$ increasing from zero, and is much better than SPC with a suitable selection of $\Lambda$. Moreover, FR-DDPC surpasses R-DDPC for large $\Lambda$, showing its ability to reduce conservatism. \textcolor{black}{Meanwhile, the value of $\Lambda$ obtained by \eqref{equation: empirically calculate Lambda} is near-optimal, which justifies the use of \eqref{equation: empirically calculate Lambda} for determining $\Lambda$. As for the theoretical bound in Theorem \ref{Theorem: bounded Lambda}, it is too large ($>10^8$ in this case) to be displayed in Fig. \ref{fig: SNR30 Lambda}.}

\section{Conclusion}
In this work, we proposed a new robust DDPC method to address the uncertainty in dynamic systems from a pessimistic perspective. Having gained insight into the inherent lack of robustness in SPC and PBR-DDPC, we constructed an uncertainty set to capture all admissible output trajectories that deviate from the output predictions by the SPC predictor and optimize the worst-case performance within the set. We established a performance guarantee for the robust DDPC method under bounded additive noise and derived a tractable convex programming reformulation with reduced-dimensional LMIs. Furthermore, integrating an affine feedback policy successfully reduced conservatism and enhanced control performance. A case study showed that a similar performance was achieved by SPC and PBR-DDPC, while the proposed method yields obvious control performance improvement. Following this work, a meaningful future step is to investigate the issue of recursive feasibility and stability.



\section*{References}
\bibliographystyle{IEEEtran}
\bibliography{ref_bib}

\end{document}